\documentclass[12pt]{article}
\usepackage{amsmath,amssymb,amsthm,eucal}
\usepackage{hyperref}
\usepackage{color}
\usepackage[normalem]{ulem}
\usepackage{enumitem}

\newcommand{\opn}{\operatorname}

\font\tenrm=cmr10

\font\cmssl=cmss10 at 12 pt

\font\bigss=cmssdc10 scaled 2300

\font\cmsslll=cmss10 at 14 pt

\renewcommand{\a}{\alpha}
\renewcommand{\b}{\beta}

\renewcommand{\d}{\delta}
\newcommand{\e}{\epsilon}

\newcommand{\g}{\gamma}

\renewcommand{\o}{\omega}

\newcommand{\x}{\xi}
\newcommand{\z}{\zeta}

\newcommand{\bC}{\mathbb{C}}
\newcommand{\bR}{\mathbb{R}}
\newcommand{\bZ}{\mathbb{Z}}

\newcommand\Sp{\mathrm{Sp}}

\newcommand{\p}{\partial}

\renewcommand{\square}{\kern1pt\vbox
               {\hrule height 0.6pt\hbox{\vrule width 0.6pt\hskip 3pt
    \vbox{\vskip 6pt}\hskip 3pt\vrule width 0.6pt}\hrule height0.6pt}
    \kern1pt}

\newcommand{\ra}{\rightarrow}

\newcommand{\ol}{\overline}

\newtheorem{Pb}{Problem}

\newtheorem{Th}{Theorem}
\newtheorem{Ex}[Th]{Example}
\newtheorem{Prop}[Th]{Proposition}

\newtheorem{Cor}[Th]{Corollary}
\newtheorem{Lem}[Th]{Lemma}
\newtheorem{Def}[Th]{Definition}

\newcommand{\bP}{\begin{Pb}\ \ }
\newcommand{\eP}{\end{Pb}}
\newcommand{\bt}{\begin{Th}\ \ }
\newcommand{\et}{\end{Th}}
\newcommand{\bp}{\begin{Prop}\ \ }
\newcommand{\ep}{\end{Prop}}
\newcommand{\bc}{\begin{Cor}\ \ }
\newcommand{\ec}{\end{Cor}}
\newcommand{\bl}{\begin{Lem}\ \ }
\newcommand{\el}{\end{Lem}}
\newcommand{\bd}{\begin{Def}\ \ }
\newcommand{\ed}{\end{Def}}

\newcommand{\pf}{\begin{proof}[{\it Proof:\ \ }]}
\newcommand{\epf}{\end{proof}}
\newcommand{\n}{\nabla}

\newcommand{\be}{\begin{equation}}
\newcommand{\ee}{\end{equation}}

\newcommand\re[1]{(\ref{#1})}
\newcommand{\arr}{\begin{array}{rlll}}
\newcommand{\ea}{\end{array}}
\newcommand{\bea}{\begin{eqnarray}}
\newcommand{\eea}{\end{eqnarray}}
\newcommand{\bean}{\begin{eqnarray*}}
\newcommand{\eean}{\end{eqnarray*}}

\catcode`@=11
\@addtoreset{equation}{section}
\catcode`@=12

\renewcommand{\z}{\zeta}
\newcommand{\zt}{\tilde{\zeta}}

\newtheorem{Rem}{Remark}  
\newcommand{\br}{\begin{Rem}\ \ }  
\newcommand{\er}{\end{Rem}}

\begin{document}
\rightline{}
\vskip 1.5 true cm
\begin{center}
{\bigss  Completeness of  projective special K\"ahler and quaternionic K\"ahler manifolds}\\[.5em]
\vskip 1.0 true cm
{\cmsslll  V.\ Cort\'es$^1$, M.\ Dyckmanns$^1$ and S.\ Suhr$^2$
} \\[10pt]
{\it\small Dedicated to Simon Salamon on the occasion of his 60th birthday}\\[10pt]
$^1${\tenrm   Department of Mathematics\\
and Center for Mathematical Physics\\
University of Hamburg\\
Bundesstra{\ss}e 55,
D-20146 Hamburg, Germany}\\
$^2${\tenrm D\'epartment de math\'ematiques et applications\\
\'Ecole normale sup\'erieure\\
45 rue d'Ulm, 75005 Paris, France\\
vicente.cortes@uni-hamburg.de, stefan.suhr@ens.fr}\\[1em]   
\vspace{2ex}
{December 27, 2016}
\end{center}
\vskip 1.0 true cm
\baselineskip=18pt
\begin{abstract}
\noindent
We prove that every projective special K\"ahler manifold with \emph{regular boundary behaviour} is complete
and defines a family of complete quaternionic K\"ahler manifolds depending
on a parameter $c\ge 0$. We also show that, irrespective of its boundary behaviour, every complete projective special K\"ahler manifold with 
\emph{cubic prepotential} gives rise to such a family. Examples include non-trivial deformations of non-compact symmetric quaternionic K\"ahler manifolds.  \\[.5em]
\end{abstract}

\tableofcontents

\section*{Introduction}

Quaternionic K\"ahler manifolds constitute a much studied class of Einstein manifolds of special holonomy \cite{B}.  
All known complete examples of positive scalar curvature are symmetric  of compact type (Wolf spaces) and it has been conjectured
that there are no more complete quaternionic K\"ahler manifolds of positive scalar curvature \cite{LS}. 
Besides the noncompact duals of the Wolf spaces, there exist also nonsymmetric complete examples of negative scalar curvature including locally 
symmetric spaces, nonsymmetric homogeneous spaces (Alekseevsky spaces) and deformations of quaternionic 
hyperbolic space \cite{L}. Our work is motivated by the desire to obtain further complete examples of 
quaternionic K\"ahler manifolds using ideas from supergravity and string theory. 

Based on general supersymmetry arguments \cite{BW} and dimensional reduction in field theory 
it has been known for a long time in the physics community that projective special K\"ahler manifolds 
(see Definition \ref{PSKDef}) are related to quaternionic K\"ahler manifolds of negative scalar curvature.
This correspondence, known as the supergravity c-map, was established by 
Ferrara and Sabharwal  \cite{FS} who explicitly associated a quaternionic K\"ahler metric with every 
projective special K\"ahler domain (see Definition \ref{PSKDDef}), cf.\ \cite{Hi} for another proof.  It was shown in \cite{CHM} that the 
supergravity c-map maps every complete projective special K\"ahler manifold to a complete
quaternionic K\"ahler manifold. 

Motivated by the fact that in the low energy limit string theory is described by supergravity, Robles Llana, Saueressig and Vandoren \cite{RSV} 
proposed a deformation of the Ferrara-Sabharwal metric (or supergravity c-map metric) depending on a real parameter.  
This deformation, know as the one-loop deformation, is interpreted as the full perturbative quantum correction (with no higher loop corrections) of supergravity when embedded into string theory. It was proven in \cite{ACDM} using  an indefinite version of the HK/QK correspondence \cite{ACM} that the one-loop deformation 
of the Ferrara-Sabharwal metric is indeed quaternionic K\"ahler on its domain of positivity. As a corollary, one obtains a new proof of the quaternionic K\"ahler property for 
the (undeformed) Ferrara-Sabharwal metric. It was also found that the completeness of the metric depends 
on the sign of the deformation parameter. In particular, it was shown that the one-loop deformation of the complex hyperbolic plane 
is complete for positive deformation parameter and incomplete for negative deformation parameter.  

The purpose of this paper is 
to give general completeness results for projective special K\"ahler manifolds and 
one-loop deformations of Ferrara-Sabharwal metrics. These results make  
it possible to construct many new explicit complete quaternionic K\"ahler manifolds of negative scalar curvature by the
supergravity c-map and its one-loop quantum correction. 

After reviewing some basic definitions and facts concerning special K\"ahler manifolds in the first section, we 
introduce the notion of regular boundary behaviour for special K\"ahler manifolds in the second section. 
The main result  of that section is that every projective special K\"ahler manifold with regular boundary behaviour is complete, see Theorem \ref{2ndThm} 
and  its Corollary \ref{1stThm}  for projective special K\"ahler domains. 

In the third section we study the one-loop deformation of Ferrara-Sabharwal metrics for nonnegative deformation parameter. We show that 
the one-loop deformation is not only defined in the case of projective special K\"ahler domains but is a globally defined one-parameter family of 
quaternionic K\"ahler metrics for every projective special K\"ahler manifold, see 
Theorem  \ref{globalFS}. Moreover, we show that the resulting quaternionic K\"ahler manifolds carry a globally defined integrable complex structure subordinate to the quaternionic structure. 

In the fourth section we prove the completeness of the  one-loop deformation for nonnegative deformation parameter under 
the assumption that the initial projective special K\"ahler manifold has either regular boundary behaviour (see Theorem \ref{DFScompleteness}) 
or is complete with cubic prepotential (see Theorem \ref{corCompleteteDefqMap2}).  The latter projective special K\"ahler manifolds are 
precisely those which can be obtained by dimensional reduction from five-dimensional supergravity \cite{DV} with complete scalar geometry \cite{CHM}. 
The corresponding construction is known as the supergravity r-map, which maps projective special real manifolds to 
projective special K\"ahler domains. 

As the simplest\footnote{The corresponding projective special real manifold is a point.} application of Theorem \ref{corCompleteteDefqMap2} (see Example \ref{remDefG22}) 
we discuss a one-parameter deformation of the metric of the noncompact symmetric space $G_2^*/SO(4)$ by locally inhomogeneous complete quaternionic K\"ahler metrics, where $G_2^*$ denotes the noncompact real form of the complex
Lie group of type $G_2$. In fact, Theorem  \ref{corCompleteteDefqMap2} implies the completeness of the one-loop deformation for all the symmetric quaternionic K\"ahler manifolds of noncompact type with exception of the quaternionic hyperbolic spaces  (which are not in the image of the supergravity c-map) and 
the spaces  $\tilde{X}(n+1)=\frac{SU(n+1,\,2)}{S[U(n+1)\times U(2)]}$. 

Similarly, applying Theorem \ref{DFScompleteness} to the complex hyperbolic space (which is a projective special K\"ahler domain with regular boundary behaviour) we obtain the completeness of the one-parameter deformation of the remaining symmetric spaces $\tilde{X}(n+1)$ ,  see Example \ref{CHEx}. 

Based on the effective necessary and sufficient completeness criterion for projective special real manifolds provided in 
 \cite[Thm.\ 2.6]{CNS}, it is easy to construct many more examples of complete projective special K\"ahler domains with cubic prepotential (see for example \cite{CDL} and
 work in progress by  J\"ungling, Lindemann and the first two authors) and corresponding one-loop deformed 
quaternionic K\"ahler manifolds by Theorem \ref{corCompleteteDefqMap2}.

\subsubsection*{Acknowledgements}
The research leading to these results has received funding 
from the the German Science Foundation (DFG) under the Research Training Group 1670 ``Mathematics inspired by String Theory" and 
from the European Research Council under the European Union's Seventh Framework 
Programme (FP/2007-2013)/ERC Grant Agreement 307062.

V.C.\ thanks the \'Ecole Normale Sup\'erieure for hospitality and support in Paris.

\section{Preliminaries}
\subsection{Conical and projective special K\"ahler manifolds}

First we recall some basic facts and definitions of special K\"ahler geometry \cite{ACD,CM}. 

\bd
A {\cmssl conical affine special K\"ahler manifold} $(M,J,g,\nabla,\xi)$ is a pseudo-K\"ahler manifold $(M,J,g)$ endowed with a flat torsion-free connection $\nabla$ and a 
vector field $\xi$ such that 
\begin{enumerate}
\item $\nabla \o =0$, where $\o=g(J.,.)$ is the K\"ahler form,
\item $d^\n J =0$, where $J$ is considered as a $1$-form with values in $TM$,
\item $\n \xi= D\xi =\mathrm{Id}$, where $D$ is the Levi-Civita connection and 
\item $g$ is positive definite on the distribution $\mathcal{D}=\mathrm{span}\{\xi,J\xi\}$ and negative definite on $\mathcal{D}^\perp$.
\end{enumerate}
\ed 

Note that the affine special K\"ahler metric $g$ has the global K\"ahler potential $f=g(\xi,\xi)$ in the sense that 
\[\frac{i}{2} \partial \overline{\partial} f=\o.\]
Furthermore the vector fields $\xi$ and $J\xi$ generate a holomorphic homothetic action of a $2$-dimensional Abelian\footnote{Note that a (real) holomorphic vector field $X$ always 
commutes with $JX$: $\mathcal{L}_X(JX)=(\mathcal{L}_X J)X=0$} Lie algebra and $J\xi$ is a Killing vector field. 

\bp
Let $(M,J,g,\nabla,\xi)$ be a conical affine special K\"ahler manifold such that the vector fields $\xi$ and $J\xi$ generate a principal $\bC^*$-action. Then the degenerate symmetric tensor field 
\be\label{g'eq}
g':=-\frac{g}{f}+\frac{\alpha^2 + (J^*\alpha)^2}{f^2},
\ee
where $\alpha := g(\xi , \cdot ) = \frac12 df$, induces a K\"ahler metric $\bar g$ on the quotient (complex) manifold $\bar M$.
\ep

\pf
It suffices to check that the kernel of $g'$ is exactly $\mathcal{D}$, the distribution tangent to the $\bC^*$-orbits, and that $g'$ is invariant under the $\bC^*$-action. 
\epf

\bd \label{PSKDef}
A {\cmssl projective special K\"ahler manifold} $(\bar{M},\bar{g})$ is a quotient as in the previous proposition with canonical projection $\pi \colon M\to \bar M$.
\ed

Notice that the projective special K\"ahler metric is related to the tensor field \re{g'eq} by $g'=\pi^* \bar g$.

\subsection{Conical and projective special K\"ahler domains}
In this section we describe an important class of special K\"ahler manifolds, the so-called special K\"ahler domains.
It is known that every special K\"ahler manifold is locally isomorphic to a special K\"ahler domain \cite{ACD}.

Let $F : M \ra \bC$ be a holomorphic function on a $\bC^*$-invariant domain $M\subset \bC^{n+1}\setminus \{ 0 \}$
such that 
\begin{enumerate}
\item[(i)] $F$ is homogeneous of degree $2$, that is $F(az)=a^2F(z)$ for all $z\in M$, $a\in \bC^*$, 
\item[(ii)] the real matrix $(N_{IJ}(z))_{I,J=0,\ldots n}$, defined by
\[ N_{IJ}(z) := 2 \mathrm{Im}\, F_{IJ}(z) = -i(F_{IJ}(z) -\ol{F_{IJ}(z)}), \]
is of signature $(1,n)$ for all $z\in M$, where $F_I := \frac{\partial F}{\partial z^I}$, $F_{IJ} := \frac{\partial^2F}{\partial z^I\partial z^J}$ etc.,  
\item[(iii)] $f(z) :=  \sum N_{IJ}(z)z^I\bar{z}^J>0$ for all $z\in M$. 
\end{enumerate}

\bd A {\cmssl conical special K\"ahler domain} $(M,g,F)$ is a $\bC^*$-invariant domain $M\subset \bC^{n+1}\setminus \{ 0 \}$ endowed with a holomorphic function 
$F$ (called {\cmssl holomorphic prepotential}) as above and with the pseudo-Riemannian metric 
\[ g = \sum N_{IJ}dz^Id\bar{z}^J.\]
\ed
Notice that $g$ has signature $(2,2n)$ and is pseudo-K\"ahler with the K\"ahler potential $f$.
A conical special K\"ahler domain becomes a conical special K\"ahler manifold if we endow it with the 
complex structure $J$ and the position vector field $\xi$  induced from the ambient space $\bC^{n+1}$. The flat connection $\n$ is induced by the standard flat connection on $\bR^{2n+2}$ via 
the immersion $M \ni (z^0,\ldots, z^n)\mapsto \mathrm{Re}(z^0,\ldots,z^n,F_0,\ldots,F_n)$.

Next we consider the domain $\bar{M} = \pi (M) \subset \bC P^n$ which is the image of $M$ under the projection 
\[ \pi : \bC^{n+1}\setminus \{ 0 \} \ra \bC P^n.\] 
The quotient manifold $\bar{M}$ inherits a (positive definite) K\"ahler metric $\bar{g}$ uniquely determined by
\be \label{pi*metric:Eq} \pi^*\bar{g} = -\frac{g}{f}+\frac{\alpha^2 + (J^*\alpha)^2}{f^2},\ee
where $\alpha := g(\xi , \cdot ) = \frac12 df$.

\bd \label{PSKDDef} A {\cmssl projective special K\"ahler domain} $(\bar{M},\bar{g})$ is the quotient $\bar{M}$  of  a conical special K\"ahler domain $M$ 
by the natural $\bC^*$-action, endowed with its canonical
K\"ahler metric $\bar{g}$. 
\ed 

Now we describe a local K\"ahler potential for the projective special K\"ahler metric $\bar{g}$
in a neighborhood of a point $p\in \bar{M}$. This yields a local K\"ahler potential $\mathcal{K}$ for projective special K\"ahler manifolds.
Let $\lambda$ be any linear function on $\bC^{n+1}$ such that $p$ lies in the affine 
chart $\{ \lambda \neq 0\}\subset \bC P^n$. The function $\frac{f}{\lambda \bar{\lambda}}$ is homogeneous 
of degree $0$ on $M \cap \{ \lambda \neq 0\}$ and therefore well defined on $\pi (M \cap \{ \lambda \neq 0\}) =
\bar{M} \cap \{ \lambda \neq 0\}$. Then 
\[ 
\mathcal{K} := -\log \left( \frac{f}{\lambda \bar{\lambda}}\right) \]
is a K\"ahler potential for the metric $\bar{g}$ on the open subset $\bar{M} \cap \{ \lambda \neq 0\}$. 
By an appropriate choice of 
linear coordinates $(z^0,\ldots ,z^n)$ on $\bC^{n+1}$ we can assume that $\lambda = z^0$.

\section{Special K\"ahler manifolds with regular boundary behaviour}

Now we consider certain compactifications of projective special K\"ahler manifolds by adding a boundary. As a first step we consider conical affine special K\"ahler manifolds
with boundary.

\bd\label{defregbd}
A {\cmssl conical affine special K\"ahler manifold with regular boundary behaviour} is a conical affine special K\"ahler manifold $(M,J,g,\n,\xi)$ which 
admits an embedding $i\colon M\to \mathcal{M}$ into a manifold with boundary $\mathcal{M}$ such that $i(M)=\text{int}\;\mathcal{M}:=\mathcal{M}
\setminus\p\mathcal{M}$ and the tensor fields $(J,g,\xi)$ smoothly extend to $\mathcal{M}$ such that, for all boundary points $p\in \partial \mathcal{M}$, 
$f(p)=0$, $df_p\neq 0$ and $g_p$ is negative semi-definite on $\mathcal{H}_p:=T_p\partial \mathcal{M}\cap J(T_p\partial \mathcal{M})$ with kernel $\mathrm{span} \{ \x_p, 
J\xi_p\}$, where $f=g(\xi,\xi)$. 
\ed

Note that for the smooth extendability of the metric $g$ it is sufficient to assume that $J$ and $f$ smoothly extend to the boundary. Indeed
this follows from the fact that $f$ is a K\"ahler potential for $g$. 

As in the case of empty boundary, we will assume that $\xi$ and $J\xi$ generate a principal $\bC^*$-action on the manifold 
$\mathcal{M}$. Then $\bar{\mathcal{M}}=\mathcal{M}/\bC^*$ is a manifold
with boundary and its interior $\bar{M}=M/\bC^*$ is a projective special K\"ahler manifold with projective special K\"ahler metric $\bar g$.  If the manifold 
$\bar{\mathcal{M}}$ with boundary
is compact, then we will call $(\bar M,\bar g)$  a  {\cmssl projective special K\"ahler manifold with regular boundary behaviour}. 

The projective special K\"ahler domains considered in Remark \ref{Rem1} below, are examples of projective special K\"ahler manifolds with regular boundary behaviour. 

\bt \label{2ndThm} Every projective special K\"ahler manifold with regular boundary behaviour is complete.
\et 

\pf Consider the underlying conical affine special K\"ahler manifold $(M,J,g,\n,\xi)$ with regular boundary behaviour. 
We first show that $g_p$ is nondegenerate for every point $p\in \p \mathcal{M}$.
By definition of regular boundary behavior we have $g|_{\mathcal{H}_p\times\mathcal{H}_p}\le 0$ with kernel $\mathrm{span} \{ \x_p, J\xi_p\}$. 
Let $\mathcal{H}_p'\subset \mathcal{H}_p$ be a complex hyperplane not containing $\x_p$. Then $g_p$ is negative definite on $\mathcal{H}_p'$. For dimensional 
reasons $\mathcal{H}_p$ is a real codimension one subspace of $T_p \partial\mathcal{M}$. Let $w$ be a vector in the complement of $\mathcal{H}_p$ in 
$T_p \partial\mathcal{M}$. By applying the Gram-Schmidt procedure we can assume that $w$ is $g_p$-orthogonal to $\mathcal{H}_p'$ in $T_p \partial\mathcal{M}$.
Then $\mathrm{span} \{ w, Jw\}$ is $g_p$-orthogonal to $\mathcal{H}_p'$ by the $J$-invariance of $g_p$. Since the real $4$-dimensional vector space $\mathrm{span} \{ \x_p, J\xi_p,w,Jw\}$ is 
$g_p$-orthogonal to $\mathcal{H}_p'$ in $T_p\mathcal{M}$ it suffices to show that $g_p$ is nondegenerate on $\mathrm{span} \{ \x_p, J\xi_p,w,Jw\}$. 
By continuity of $df$ and $\x$ we know that 
\[2g_p(\x_p,.)=df_p.\]
Since $Jw\notin T_p\partial \mathcal{M}$ and $w\in T_p\partial\mathcal{M}$ we have 
\[0\neq df_p(Jw)=2g_p(\x_p,Jw)=-2g_p(J\x_p,w)\text{ and }0=df_p(w)=2g_p(\x_p,w)=2g_p(J\x_p,Jw).\]
Now by considering the representing matrix of $g_p$ on $\mathrm{span} \{ \x_p, J\xi_p,w,Jw\}$ and using that $g_p$ vanishes on $\mathrm{span} \{ \x_p, J\xi_p\}$ 
we see that $g_p$ is nondegenerate. 
This proves that $g_p$ is nondegenerate and, therefore, of signature $(2,2n)$ by continuity. 

Let $\g : I \ra \bar{M}$, $I=[0, b)$, $0<b\le \infty$, be a curve which is not contained in any compact subset of $\bar{M}$. We will show that $\gamma$ has infinite length under the assumption of regular boundary behaviour. 
Call a point $p\in \bar{\mathcal{M}}$ an accumulation point of $\g$ if there exists a sequence $t_i\in I$ such that $\lim t_i=b$ and $\lim \g(t_i)=p$. By our assumption, $\g$ has at least 
one accumulation point $\bar{p}_0$ on the boundary.  We distinguish two cases:

1$^{st}$ case: $\gamma$ has exactly one accumulation point $\overline{p}_0$ which necessarily lies on the boundary.
Under this hypothesis, for every neighborhood of $\bar{p}_0$ we can find $a\in I$ such that $\g([a,b))$ is fully contained in that neighborhood.

Choose a point $p_0\in \pi^{-1}(\bar{p}_0)\subseteq \p \mathcal{M}$. Since the signature of $g_{p_0}$ is $(2,2n)$, there exists a complex hyperplane $E \subset 
T_{p_0}\mathcal{M}$ on which 
$-g$ is positive definite. Let $M'$ denote a complex hypersurface through $p_0$ tangent to $E$ such that $-g|_{TM'\times TM'}$ is positive definite.

The pullback of the projective special K\"ahler metric can be estimated on $N=\mathrm{int}(M')$ as follows
\be \label{est2:Eq} (\pi^*\bar{g})|_N  = -\left.\frac{g}{f}\right|_N+ \left. \frac{\alpha^2 + (J^*\alpha)^2}{f^2}\right|_N  \ge 
\left. \frac{\alpha^2}{f^2}\right|_N = \frac{df^2}{4f^2}.\ee
Now we show how this implies that $\g$ has infinite length. We can assume by shifting the initial point of the interval $I$ that $\g$ is fully contained 
in $\pi(N)\subset \bar{M}$. Let $\g_N : I \ra N$ be the curve which 
projects to $\g$ under $\pi|_N$.  Then there exists 
a sequence $t_i \in [0, b)$ such that $f(\g_N(t_i)) \ra 0$ and $\g_N([0,t_i])\subset \g_N(I) \subset N$.
In view of \re{est2:Eq}, we have
\begin{align*}
L(\g ) &\ge L\big(\g|_{[0,t_i]}\big)
=  L^{ \pi^*\bar{g}} (\g_N|_{[0,t_i]})
\ge \frac12 \int_{0}^{t_i}\left| \frac{d}{dt}\log f \circ \g_N\right| dt\\
&\ge -\frac12 \int_{0}^{t_i}\frac{d}{dt}\log f \circ \g_N \;dt
= \frac12 \big(\log f(\g_N(0)) - \log f(\g_N(t_i))\big) \ra \infty .
\end{align*}
This shows that $\g$ has infinite length. 

2$^{nd}$ case: $\g$ has at least two accumulation points. Let $\overline{p}_0\neq 
\overline{p}_1$ be such accumulation points. We know that at least one accumulation point, e.g.\ 
$\overline{p}_0$, lies in the boundary. Under the assumption that there exists a second accumulation 
point, we now show that the second accumulation point can be taken arbitrarily near to 
$\overline{p}_0$. In other words, we claim that for every given neighborhood $U$ of $\bar{p}_0$ there exists an accumulation 
point $\bar{p}_2\in U\setminus \{\bar{p}_0\}$. Indeed let us denote by $B^{aux}_r(\bar{p}_0)$ the ball 
of radius $r>0$ centered at $\bar{p}_0$ with respect to an auxiliary Riemannian metric on
$\bar{\mathcal{M}}$. 
Choose $r>0$ such that $\overline{B_r^{aux}(\bar{p}_0)}\subset U$. If 
$\bar{p}_1\in U$ there is nothing to prove. If $\bar{p}_1\notin U$ choose 
sequences $s_i < t_i < s_{i+1}$ such that $\lim_{i\ra \infty}  \g (s_i)=\bar{p}_0$ and $\lim_{i\ra 
\infty}  \g (t_i)=\bar{p}_1$. 
We can assume that $\g(s_i)\in  B_{r/2}^{aux}(\bar{p}_0)$ and $\g(t_i)\notin \overline{B_{r}^{aux}(\bar{p}_0)}$ 
for all $i$. Then there exists a
sequence $u_i \in (s_i,t_i)$ with $\g(u_i) \in B_r^{aux}(\bar{p}_0)
\setminus B_{r/2}^{aux}(\bar{p}_0)$. 
The sequence $\g(u_i)$ has an accumulation point $\bar{p}_2\in 
\overline{B_r^{aux}(\bar{p}_0)}\subset U$. We will continue to denote this accumulation point 
arbitrarily close to $\overline{p}_0$ by $\bar{p}_1$.

If $\overline{p}_1 \in \overline{M}$ it is easy to see that $\g$ has infinite length. In fact 
consider a geodesically convex ball $B_\d (\bar{p}_1)$ of radius $\d >0$ centered at $\bar{p}_1$ with
respect to $\bar{g}$. We take $\d$ sufficiently small such that $B_\d (\bar{p}_1)$ is relatively compact
in $\bar{M}$. Since the curve $\g$ intersects the ball  $B_{\d /2} (\bar{p}_1)$ an arbitrarily large 
number of times $k$, the length of $\g$ is larger or equal than $k\d \ra \infty$.

Thus we can assume that $\bar{p}_1$ lies in the boundary as well. By restricting $U$ we can assume
that $U$ is in the image of a complex hypersurface $M'\subset \overline{\mathcal{M}}$ as above. 
We can further assume that $f\le \e$ on $M'$. 
Since $g'=\pi^*\bar g$ is given by \re{g'eq} the Riemannian metric $\pi^*\bar{g}|_N$ on $N=\mathrm{int}(M')$ is bounded from below by the 
Riemannian metric 
\be \label{gF2:Eq} -\left. \frac{g}{f}\right|_N \ge -\frac{1}{\e}g|_N.\ee
Let us denote by $B'_r(p)$ the ball centered at $p\in M'$ of radius $r>0$ with respect to the Riemannian metric $-g|_{M'}$ on $M'$. 
We choose $\d >0$ such that $B'_\d(p_0)$ is relatively compact in $M'$. Then 
every curve in $B'_{\d}(p_0)$ from $B'_{\d/2}(p_0)\subset M'$ $(p_0:=(\pi|_{M'})^{-1}(\bar{p}_0))$ which leaves $B'_{\d}(p_0)$ has length 
with respect to  $-g|_{M'}$  bounded from below by some positive constant $c$ (in fact
$c=\d/2$). Since we can assume that $p_1:=(\pi|_{M'})^{-1}(\bar{p}_1)$ is arbitrarily close to $p_0$ we can assume that 
$p_1 \in B'_\d(p_0)$ and there exist disjoint balls $B'_{\d'}(p_0), B'_{\d'}(p_1)\subset B'_\d(p_0)$ which have distance with respect to $-g|_{M'}$ bounded from below
by some positive constant. By reducing the above constant $c$, if necessary, we can assume that this constant is again $c$.   
Then we can conclude that every curve which connects a point in 
$B'_{\d'}(p_0)$ with a point in $B'_{\d'}(p_1)$ has length with respect to $-g|_{M'}$  bounded from below by $c$. 
Since $p_0$ and $p_1$ are accumulation points of $\g$ either $\g$ leaves the set 
$\pi (N)$ infinitely often, in which case $\g$ has infinite length, or $\g$ stays eventually inside
$\pi (N)$, in which case it can be eventually identified with a curve $\g_N$ in $N$ by the projection
$\pi|_{N}$. Since $\bar{p}_0$ and $\bar{p}_1$ are accumulation points of $\g_N$ there exists an infinite number
of arcs of $\g_N$ in $N$ connecting $B'_{\d'}(p_0)$ with  $B'_{\d'}(p_1)$. Again the 
length is infinite. In both cases we used the estimate \re{gF2:Eq} together with the lower bound $c$ on the length 
of arcs with respect to $-g|_{N}$. 
\epf

\begin{Rem}\label{Rem1} 
In the case of conical affine special K\"ahler domains the description of regular boundary behaviour simplifies as follows. 
Let $(\bar{M},\bar{g})$ be a projective special K\"ahler domain with underlying conical special K\"ahler domain 
$(M,g,F)$. Suppose that the affine K\"ahler potential $f$ extends to a smooth function (denoted again by $f$) 
on some neighborhood of $\mathrm{cl}(M)\setminus \{ 0\}$, where $\mathrm{cl}(M)$ denotes the 
closure of $M$, such that $f(p)=0$, $df_p\neq 0$, and that $g_p$ is negative semi-definite on $T_p\partial M\cap J(T_p\partial M)$ 
with kernel $\bC \xi_p=\bC p$ for all boundary points $p\in \partial M \setminus \{ 0\}$. Then $(M,g,F)$ 
is an example of a conical affine special K\"ahler manifold with regular boundary behaviour and 
$(\bar{M},\bar{g})$ an example of a projective special K\"ahler manifold with regular boundary behaviour.
\end{Rem}

The following result is an immediate consequence of Theorem \ref{2ndThm}.

\bc \label{1stThm} Under the above assumptions on the boundary behaviour of the affine K\"ahler potential $f$ in Remark \ref{Rem1},  the Riemannian manifold $(\bar{M},\bar{g})$
is complete. 
\ec

\section{One-loop deformed Ferrara-Sabharwal metric}

In this section we will recall the definition of the one-loop (quantum) deformation of the Ferrara-Sabharwal metric which is a one-parameter family of 
quaternionic K\"ahler metrics associated with a projective special K\"ahler domain \cite{RSV,ACDM}. The fact that the metric is quaternionic K\"ahler 
was proven in \cite{ACDM} with the help of an indefinite version of Haydys' HK/QK correspondence \cite{Haydys} developed in \cite{ACM}. This implies 
that the reduced scalar curvature $\nu=\frac{\mathrm{scal}}{4m(m+2)}$ is negative and more precisely given by $\nu =-2$ with the present normalizations.
Here $m$ is the quaternionic dimension of the quaternionic K\"ahler manifold. In the special case of the (undeformed) Ferrara-Sabharwal metric the quaternionic 
K\"ahler property was obtained by different methods in \cite{FS,Hi}.

Every projective special K\"ahler manifold admits 
a covering by projective special K\"ahler domains and we will show that the one-loop deformed Ferrara-Sabharwal metrics associated with the domains
can be consistently glued to a globally defined (quaternionic K\"ahler; to be shown) metric. This generalizes the result that the Ferrara-Sabharwal 
metric, which was originally defined for special K\"ahler domains \cite{FS}, is globally defined for every projective special K\"ahler manifold \cite{CHM}. 
We will also show that the above quantum deformed quaternionic K\"ahler manifolds admit a globally defined integrable complex structure
$J_1$ subordinate to the quaternionic structure, generalizing results of \cite{CLST} for the Ferrara-Sabharwal metric.

\subsection{The supergravity c-map}
\label{sugraSect}
Let $(\bar{M},\,\bar{g})$ be a projective special K\"ahler domain of complex dimension $n$. 
The \textbf{supergravity c-map} \cite{FS} associates with $(\bar{M},\,\bar{g})$ 
a quaternionic K\"ahler manifold $(\bar{N},\,g_{\bar{N}})$ of dimension $4n+4$. 
Following the conventions of \cite{CHM}, we have $\bar{N}=\bar{M}\times \mathbb{R}^{>0}\times \mathbb{R}^{2n+3}$ and
\begin{eqnarray*} \label{e:2} g_{\bar{N}} &=& \bar{g} + g_G,\\\nonumber 
g_G&=& \frac{1}{4\rho^2}d\rho^2 + \frac{1}{4\rho^2}\left(d\tilde{\phi}
+ \sum \left(\zeta^Id\tilde{\zeta}_I-\tilde{\zeta}_Id\zeta^I\right) \right)^2 
+\frac{1}{2\rho}\sum \mathcal{I}_{IJ}(m) d\zeta^Id\zeta^J\\ 
&&+ \frac{1}{2\rho}\sum \mathcal{I}^{IJ}(m)(d\tilde{\zeta}_I + 
\mathcal{R}_{IK}(m)d\zeta^K)(d\tilde{\zeta}_J + 
\mathcal{R}_{JL}(m)d\zeta^L),
\end{eqnarray*}
where $(\rho,\,\tilde{\phi},\,\tilde{\z}_I,\, \z^I)$, $I=0,\,1,\,\ldots,\,n$, 
are standard coordinates on $\bR^{>0}\times \bR^{2n+3}$. 
The real-valued matrices $\mathcal{I}(m):=(\mathcal{I}_{IJ}(m))$ and $\mathcal{R}(m):=(\mathcal{R}_{IJ}(m))$
depend only on $m\in \bar{M}$ and $\mathcal{I}(m)$ is invertible with
the inverse $\mathcal{I}^{-1}(m)=:(\mathcal{I}^{IJ}(m))$. More precisely,
 \begin{equation} 
\label{FRI}
{\cal N}_{IJ} := 
\mathcal{R}_{IJ} + i\mathcal{I}_{IJ} := 
\bar{F}_{IJ} + 
i \frac{\sum_K N_{IK}z^K\sum_L N_{JL}z^L}{\sum_{IJ}N_{IJ}z^Iz^J} ,\quad 
N_{IJ} := 2 \mathrm{Im} F_{IJ} ,\end{equation}
where $F$ is the holomorphic prepotential with respect
to some system of special holomorphic coordinates $(z^I)$ on the 
underlying conical special K\"ahler domain $M\to \bar{M}$. 
Notice that the expressions are homogeneous of degree zero and, hence, 
well-defined functions on $\bar{M}$. It is shown in \cite[Cor.\ 5]{CHM} 
that the matrix $\mathcal{I}(m)$ is positive definite
and hence invertible and that the metric $g_{\bar{N}}$ does not
depend on the choice of special coordinates \cite[Thm.\ 9]{CHM}. 
It is also shown that $(\bar{N},\,g_{\bar{N}})$ is complete if 
and only if $(\bar{M},\,\bar{g})$ is complete \cite[Thm.\ 5]{CHM}. 
\noindent
Using $(p_a)_{a=1,\,\ldots,\, 2n+2}:=(\zt_I,\,\z^J)_{IJ=0,\ldots,n}$ and the positive definite matrix \cite{CHM}
$$(\hat{H}^{ab}):=\begin{pmatrix}\mathcal{I}^{-1} & \mathcal{I}^{-1}\mathcal{R} \\ \mathcal{R}\mathcal{I}^{-1} & \mathcal{I}+\mathcal{R}\mathcal{I}^{-1}\mathcal{R}\end{pmatrix},$$
we can combine the last two terms of $g_G$ into $\frac{1}{2\rho}\sum dp_a \hat{H}^{ab} dp_b$, i.e.\ the quaternionic K\"ahler metric is given by
\begin{equation}\label{FSmetric} g_{FS}:=g_{\bar{N}}=\bar{g}+\frac{1}{4\rho^2}d\rho^2 + \frac{1}{4\rho^2}\left(d\tilde{\phi}
+ \sum \left(\zeta^Id\tilde{\zeta}_I-\tilde{\zeta}_Id\zeta^I\right) \right)^2+ \frac{1}{2\rho}\sum dp_a \hat{H}^{ab} dp_b.\end{equation}

\noindent
This metric is known as the \textbf{Ferrara-Sabharwal} metric.

\subsection{The one-loop deformation}\label{secHKQKCMapNew}

Now we consider a family of metrics $g^c_{FS}$ depending on a real parameter $c$ such that $g^0_{FS}=g_{FS}$. To define this family we assume for the moment that $z^0\neq 0$
on the conical affine special K\"ahler domain $M\subset \bC^{n+1}$. Under this assumption we can  consider the projective special K\"ahler domain as a subset $\bar M\subset 
\bC^n\subset \bC P^n$.

\bd\label{defDeformedFSMetric}
For any $c\in\mathbb{R}$, the metric 
\begin{eqnarray}g^c_{FS}&=&\frac{\rho+c}{\rho}\, \bar{g} +\frac{1}{4\rho^2}\frac{\rho+2c}{\rho+c}d\rho^2+\frac{1}{4\rho^2}\frac{\rho+c}{\rho+2c}(d\tilde{\phi}+\sum_{I=0}^n(\z^Id\zt_I-\zt_Id\z^I)+cd^c\mathcal{K})^2\nonumber \\
&&+\frac{1}{2\rho}\sum_{a,\,b=1}^{2n+2} dp_a\hat{H}^{ab}dp_b+\frac{2c}{\rho^2}e^{\mathcal{K}}\left|\sum_{I=0}^n (X^Id\zt_I+F_I(X)d\z^I)\right|^2\label{DefFSmetric2}
\end{eqnarray}
is defined on the domains
\begin{align}
N'_{(4n+4,\,0)}&:=\{\rho>-2c,~\rho>0\}\subset \bar{N},\nonumber\\
N'_{(4n,\,4)}&:=\{-c<\rho<-2c\}\subset\bar{N},\nonumber\\ N'_{(4,\,4n)}&:=\bar{M}\times \{-c<\rho<0\}\times\mathbb{R}^{2n+3}\subset \bar{M}\times \mathbb{R}^{<0}\times \mathbb{R}^{2n+3}\label{eqFSDomain}
\end{align}
for any projective special K\"ahler domain $\bar{M}$ defined by a holomorphic function $F$ on the underlying conical affine special K\"ahler domain 
$M$, where $\bar{N}=\bar{M}\times\mathbb{R}^{>0}\times \mathbb{R}^{2n+3}$, $(X^\mu)_{\mu=1,\,\ldots,\,n}$ are standard inhomogeneous holomorphic coordinates on $\bar{M}
\subset \bC^n$, $X^0:=1$, the real coordinate $\rho$ corresponds to the second factor,  $(\tilde{\phi},\,\tilde{\zeta}_I,\,\zeta^I)_{I=0,\,\ldots,\,n}$ are standard real coordinates on $\mathbb{R}^{2n+3}$, and $\mathcal{K}:=-\log \sum_{I,\,J=0}^n X^IN_{IJ}(X)\bar{X}^J$ is the K\"ahler potential for $\bar{g}$.
The metric $g^c_{FS}$ is called the \textbf{one-loop deformed Ferrara-Sabharwal metric}.
\ed

\bp\label{propDefFSPosParamIsometric}
Let $\bar{M}\subset \bC^n\subset \bC P^n$ be a projective special K\"ahler domain and $g^c_{FS}$, $g^{c'}_{FS}$ one-loop deformed Ferrara-Sabharwal metrics for positive deformation parameters $c,\,c'\in\mathbb{R}^{>0}$ defined on $\bar{N}=N'_{(4n+4,\,0)}$. Then $(\bar{N},\,g^{c}_{FS})$ and $(\bar{N},\,g^{c'}_{FS})$ are isometric.
\ep
\pf
Any $e^{\lambda}\in\mathbb{R}^{>0}$ acts diffeomorphically on $\bar{N}=\bar{M}\times\mathbb{R}^{>0}\times \mathbb{R}^{2n+3}$ as follows:
\[\bar{N}\to\bar{N},\quad (m,\,\rho,\,\tilde{\phi},\,\tilde{\zeta}_I,\,\zeta^I)_{I=0,\,\ldots,\,n}\mapsto (m,\,e^{\lambda}\rho,\,e^{\lambda}\tilde{\phi},\,e^{\lambda/2}\tilde{\zeta}_I,\,e^{\lambda/2}\zeta^I)_{I=0,\,\ldots,\,n}.\]
Under this action, $g^c_{FS}\mapsto g^{e^{-\lambda}c}_{FS}$. Choosing $e^{\lambda}=c/c'$, this shows that \linebreak $(\bar{N},\,g^{c}_{FS})$ and $(\bar{N},\,g^{c'}_{FS})$ are isometric.
\epf

\subsection{Globalization of the one-loop deformed metric} 

Let $(\bar M, \bar g)$ be a projective special K\"ahler manifold with underlying conical affine special K\"ahler manifold $(M,J,g,\nabla, \xi)$. Consider a covering of 
$\bar M$ by open subsets $\bar{M}_\alpha$ isomorphic to projective special K\"ahler domains. Over the preimage $M_\alpha:=\pi^{-1}(\bar{M}_\a)$ we have a system
of so-called conical affine special coordinates $(z^I)_{0\le I\le n}$ which correspond to the natural coordinates in the underlying conical affine special K\"ahler domain
equipped with the holomorphic prepotential $F$. Notice that the map $\phi_\a\colon M_\a \rightarrow \bC^{2n+2}$, $p\mapsto (z^I, F_{I})|_p$, where $F_I$ denotes the 
$I$-th partial derivative at the point $z=(z^0,\ldots, z^n)$, is a conical nondegenerate Lagrangian immersion in the sense of \cite{CM}. Further note that the coordinates as well as the prepotential depend on $\a$. To indicate this dependence we will write $z^I_{\a}$, $F_{\a}$ etc. Since any pair of conical 
nondegenerate Lagrangian immersions is related by a real linear symplectic transformation \cite{ACD, CM} there exists an element 
$$\mathcal{O}=\mathcal{O}_{\beta,\a}=\left(\begin{matrix}
A & B \\ C & D
\end{matrix}\right)\in \Sp(\bR^{2n+2})$$ 
such that $\phi_\beta =\mathcal{O}\circ \phi_\a$ on $M_\a\cap M_\beta$. 

Define $\bar{N}_\a:= \bar{M}_\a\times \bR^{>0}\times S^1_c\times \bR^{2n+2}$ and $N_\a := M_\a \times \bR^{>0}\times S^1_c\times \bR^{2n+2}$,  where 
$\bR^{>0}\times\bR^{2n+3}$ is endowed with the standard coordinate system $(\rho,\tilde\phi,
\zt_I,\zeta^J)=(\rho_\a,\tilde\phi_\a,\zt_{I,\a},\zeta^J_\a)=:(\rho_\a,\tilde\phi_\a,v_\a)$ and $S^1_c:=\bR/2\pi c\bZ$. Notice that $S^1_c$ can be canonically identified with $S^1=\bR/2\pi \bZ$ by $[x]\mapsto [cx]$ if $c\neq0$ and that $S^1_0=\bR$.

Next we define an equivalence relation on the disjoint union of the $\bar{N}_\a$ (and similarly on the disjoint 
union of the $N_\a$)
\begin{eqnarray*}
&(m_\a,\rho_\a,\tilde\phi_\a,v_\a)\sim (m_\beta,\rho_\beta,\tilde\phi_\beta,v_\beta)\\
&:\Leftrightarrow m_\a=m_\beta,\; \rho_\a=\rho_\beta, \;
\tilde\phi_\beta=\tilde\phi_\a - ic\log\left(\frac{z^0_\a\bar{z}^0_\beta}{z^0_\beta \bar{z}^0_\a}\right),\; v_\beta =(\mathcal{O}^t_{\b,\a})^{-1} v_\a.
\end{eqnarray*}
\bp\label{QKGL}
The quotient $\bar N:=\cup_\a \bar{N}_\a/\sim$ is a smooth manifold of real dimension $4n+4$ fibering over the projective special K\"ahler manifold $\bar M$ as a 
bundle of  flat symplectic manifolds modeled on the quotient of a symplectic vector space $\bR^{2n+2}$ by a cyclic group of translations (the cyclic group is trivial for $c=0$). By $\pi$ we denote the induced natural projection $\bar{N}\to\bar{M}$. Similarly, the quotient  $N:=\cup_\a N_\a/\sim$ is a bundle over the conical affine 
special K\"ahler manifold $M$ with flat symplectic fibers.
\ep

\pf
It is clear that $\bar N$ is a fibre bundle with standard fibre $\bR^{>0}\times S^1_c\times  \bR^{2n+2}$. By taking the logarithm of $\rho$ one can identify the standard
fibre with  the quotient $\bR\times S^1_c\times \bR^{2n+2}$ of $\bR^{2n+4}$ by the group of translations $2\pi c \bZ$ acting on the second coordinate. Since the transition functions take values in the group of affine symplectic transformations of $\bR\times S^1_c\times \bR^{2n+2}$, the  fibers of the resulting bundle naturally  carry a flat symplectic structure. In fact, the linear part of the transition functions takes values in the subgroup $\{\mathrm{Id}_{\bR^2}\}\times \Sp(\bR^{2n+2})
\subset \Sp(\bR^{2n+4})$. 
\epf

To avoid a parameter-dependence of the domain of definition of the metric we will assume from now on for simplicity that the one-loop parameter $c>0$.

\bt\label{globalFS}
The quaternionic K\"ahler metrics $g_{FS,\a}^c$, $c>0$, given by \eqref{DefFSmetric2} on each coordinate domain $\bar{N}_\a$ of $\bar N$ using the coordinates 
$(X^\mu,\rho,\tilde\phi,\zt_{I},\zeta^J)=(X^\mu_\a,\rho_\a,\tilde\phi_\a,\zt_{I,\a},\zeta^J_\a)$ induce a well-defined quaternionic K\"ahler metric $g_{FS}^c$ on $\bar N$. 
Furthermore there exists a globally defined integrable complex structure $J_1$ subordinate to the parallel skew-symmetric quaternionic structure $Q$ of $(\bar N, g_{FS}^c)$.
\et

\pf
First we show that the quaternionic K\"ahler metrics defined on the domains $\bar{N}_\a$ are consistent. 
The terms $\frac{\rho+c}{\rho}\bar{g}$ and $\frac{1}{4\rho^2}\frac{\rho+2c}{\rho+c}d\rho^2$ in \eqref{DefFSmetric2} are 
manisfestly coordinate independent, since the transition functions do not act on $\rho$. The one-form $\eta_{can}:=\sum_{I=0}^n(\z^Id\zt_I-\zt_Id\z^I)$ is obviously invariant under 
linear symplectic transformations and therefore also coordinate independent. The invariance of the term $\sum_{a,\,b=1}^{2n+2} dp_a\hat{H}^{ab}dp_b$ was shown in 
\cite[Lemma 4]{CHM}. Next we show the invariance of $d\tilde{\phi}+cd^c\mathcal{K}$. Since 
\[\sum_{I,J} X^I N_{IJ} \bar{X}^J=\frac{f}{z^0\bar{z}^0},\]
where $f=g(\xi,\xi)=\sum_{I,J} z^I N_{IJ} \bar{z}^J$ is coordinate independent (but defined on $N$, not on $\bar N$), we see that 
\[cd^c\mathcal{K}_\b-cd^c\mathcal{K}_\a=cd^c\log \left(\frac{z^0_\b\bar{z}^0_\b}{z^0_\a\bar{z}^0_\a}\right)
=icd\log\left(\frac{z^0_\a\bar{z}^0_\beta}{z^0_\beta \bar{z}^0_\a}\right),\] 
where we have used that $d^c=-J^*d$ on functions. By the transition rule for $\tilde\phi$ we have
\[d\tilde\phi_\b= d\tilde\phi_\a -icd\log\left(\frac{z^0_\a\bar{z}^0_\beta}{z^0_\beta \bar{z}^0_\a}\right).\]
This shows the invariance of $d\tilde{\phi}+cd^c\mathcal{K}$. 

Finally we show the invariance of $e^{\mathcal{K}}\left|\sum_{I=0}^n (X^Id\zt_I+F_I(X)d\z^I)\right|^2$. By rewriting this as 
\begin{eqnarray*}
\frac{1}{\sum X^I N_{IJ}(X) \bar{X}^J} \left|\sum_{I=0}^n (X^Id\zt_I+F_I(X)d\z^I)\right|^2&= \frac{z^0\bar{z}^0}{f}\left|\sum_{I=0}^n \frac{z^I}{z^0}d\zt_I
+F_I(\frac{z}{z^0})d\z^I\right|^2\\
&=\frac{1}{f}\left|\sum\nolimits_{I=0}^n z^Id\zt_I +F_I(z)d\z^I\right|^2
\end{eqnarray*}
we see that the term is coordinate independent. In fact, the sum $\sum\nolimits_{I=0}^n z^Id\zt_I +F_I(z)d\z^I$ is obtained from the natural pairing between 
$\bC^{2n+2}$ and $(\bC^{2n+2})^*\supset (\bR^{2n+2})^*$ which is, in particular, invariant under linear symplectic transformations.
Summarizing we have shown that the metric $g_{FS}^c$ is well defined on $\bar N$.

Since $g_{FS}^c$ is quaternionic K\"ahler (of negative scalar curvature) on each of the domains $\bar{N}_\a$ it follows that $g_{FS}^c$ is a quaternionic K\"ahler metric. In fact, 
the locally defined parallel skew-symmetric quaternionic structures on the domains $\bar{N}_\a$ are uniquely determined by the Lie algebra of the holonomy group of 
$g_{FS}^c|_{\bar{N}_\a}$ and therefore extend to a globally defined quaternionic structure $Q$. It can be also checked by direct calculations (see below) that the locally defined quaternionic 
structures $Q_\a$ on $\bar{N}_\a$ are consistent. In fact, the description of the quaternionic K\"ahler structure on $\bar{N}_\a$ in terms of the HK/QK-correspondence 
\cite{ACDM} yields an almost hypercomplex structure $(J_1,J_2,J_3)$ on $\bar{N}_\a$ which defines the quaternionic structure $Q_\a$. The structure is defined by the 
three K\"ahler forms $\o_i= g_{FS}^c(J_i\cdot ,\cdot)$, $i=1,2,3$. These are given by 
\[\o_i=-d\theta_i+2\theta_j\wedge \theta_k,\]
where $(i,j,k)$ is a cyclic permutation of $\{1,2,3\}$ and the one-forms $\theta_i$ on $\bar{N}_\a$ are defined by 
\begin{eqnarray*}
\theta_1&=&-\frac{1}{4\rho}\big(d\tilde{\phi}+(\rho+c)d^c\mathcal{K}-\eta_{can}\big)\\
\theta_2+i\theta_3&=&i\frac{\sqrt{\rho+c\,}\,}{\rho}e^{\mathcal{K}/2} \sum_{I=0}^n X^I A_I,\quad A_I:=d\zt_I+\sum_J F_{IJ} d\zeta^J.
\end{eqnarray*}

Next we prove that $Q$ admits a global section $J_1$ by showing that the K\"ahler form $\o_1$ is invariantly defined, i.e.\ coordinate independent. 
First we remark that 
$\theta_1$ can be decomposed as 
\[\theta_1= -\frac{1}{4\rho}\big(d\tilde{\phi}+cd^c\mathcal{K}-\eta_{can}\big) -\frac{1}{4}d^c\mathcal{K},\]
where the first was already shown to be invariant. Using that $\mathcal{K}= -\log\frac{f}{(r^0)^2}$, where $z^0=r^0e^{i\varphi^0}$, the second term can be decomposed as
\[-\frac{1}{4}d^c\mathcal{K}= \frac{1}{4}d^c\log f +\frac{1}{2} d^c\log r^0=\frac{1}{4}d^c\log f -\frac{1}{2}d\varphi^0.\]
Since the first term on the right-hand side is invariant we see that 
\[\theta_1=\theta_1^{inv}-\frac{1}{2}d\varphi^0,\]
where $\theta_1^{inv}$ is coordinate independent. This implies that $d\theta_1$ is invariant. Now we observe that 
\[\sum z^I A_I\]
is invariant (defined on $N$). This follows from 
\[\sum z^I A_I= \sum z^Id\zt_I +F_I(z)d\z^I,\]
where the right-hand side was already observed to be invariant. As a consequence, the two-form 
\[\theta_2\wedge \theta_3= -\frac{1}{2i} (\theta_2+i\theta_3)\wedge (\theta_2-i\theta_3)\]
is also invariant, since 
\[\theta_2+i\theta_3= \frac{i}{\rho} \left( \frac{\rho +c}{f}\right)^\frac12 e^{-i\varphi^0}\sum z^I A_I,\]
which implies that $e^{i\varphi^0}(\theta_2+i\theta_3)$ is a well defined one-form on $N$.

Combining these results we have shown that $\o_1=-d\theta_1+2 \theta_2\wedge \theta_3$ is invariant. By similar calculations it is easy to show that a conformal multiple 
$e^{i\varphi^0}\o$ of the $(2,0)$-form
\[\o=\o_2+i\o_3\]
with respect to $J_1$ is invariantly defined on $N$ (and horizontal with respect to the projection $N\ra \bar{N}$ induced by $M\ra \bar{M}$). This implies that the complex plane spanned by $\o$ and $\bar\o$ is invariantly defined on $\bar N$ and therefore the real plane spanned
by $\o_2$ and $\o_3$. This reproves the fact that the quaternionic structure is well-defined.

Now we prove the integrability of $J_1$. It is sufficent to check this on $\bar{N}_\a$. In the case $c=0$ this was previously shown in \cite{CLST}. With the definition of $\o_1$ above
we compute
\begin{align}
\omega_1&=\frac{1}{4\rho}\left(d\rho\wedge d^c\mathcal{K}+(\rho+c)\,dd^c\mathcal{K}-2\sum_{I=0}^n d\zt_I\wedge d\z^I\right)+\frac{1}{\rho} d\rho\wedge \theta_1\nonumber\\
&\quad\quad +\frac{\rho+c}{\rho^2}e^{\mathcal{K}}i(\sum_I X^IA_I)\wedge(\sum_J \bar{X}^J\bar{A}_J)\nonumber\\
&=\frac{\rho+c}{\rho}\frac{1}{4}dd^c\mathcal{K}+\frac{i}{2}\frac{1}{4\rho^2}\frac{\rho+c}{\rho+2c}\tau\wedge\bar{\tau}-\frac{i}{2}\frac{1}{\rho}\sum_{I,\,J=0}^nN^{IJ}A_I\wedge\bar{A}_J\nonumber\\
&\quad\quad +\frac{i}{2}\frac{2\rho+2c}{\rho^2}e^{\mathcal{K}}(\sum_I X^IA_I)\wedge(\sum_J \bar{X}^J\bar{A}_J)\label{eqDefFSOmega1},
\end{align}
where
\[ \tau:= d\tilde{\phi}+\eta_{can}+cd^c\mathcal{K}+i\frac{\rho+2c}{\rho+c}d\rho\]
and we used that 
\[\sum_{I,\,J=0}^n iN^{IJ}A_I\wedge\bar{A}_J=\sum_{I,\,J,\ K=0}^n iN^{IJ}(F_{IK}-\bar{F}_{IK})d\zeta^K\wedge \zt_J
=\sum_{I=0}^n d\tilde{\zeta}_I\wedge d\zeta^I.\]

Together with the expression
\[ g^c_{FS}=\frac{\rho+c}{\rho} \bar{g}+\frac{1}{4\rho^2}\frac{\rho+c}{\rho+2c}|\tau|^2-\frac{1}{\rho}\sum_{I,\,J=0}^n N^{IJ}A_I\,\bar{A}_J+\frac{2\rho+2c}{\rho^2}e^{\mathcal{K}}\big|\sum_{I=0}^nX^IA_I\big|^2\]
 for the deformed Ferrara-Sabharwal metric, which can be proven using \cite[Lemma 3]{ACDM}, \eqref{eqDefFSOmega1} shows that 
\[(\tau,\,dX^\mu,\, A_I)_{I=0,\,\ldots,\,n}^{\mu=1,\,\ldots,\,n}\] 
is a coframe of holomorphic one-forms with respect to $J_1$. This can be linearly combined into the coframe
\begin{align*}\big(\tau+2ic\partial \mathcal{K}-&2\sum_{I=0}^n\zeta^IA_I-\sum_{I,\,J,\,K=0}^n\zeta^I F_{IJK}(X)\zeta^JdX^K,\\&\qquad\qquad~dX^\mu,
~\frac{1}{2}(A_I-\sum_{J,\,K=0}^nF_{IJK}(X)\zeta^JdX^K)\big)\end{align*}
of closed holomorphic one-forms which corresponds to the $J_1$-holomorphic coordinate system
\[ (\chi,~ X^\mu,~ w_I=\frac{1}{2}(\tilde{\zeta}_I+\sum_{J=0}^n F_{IJ}(X)\zeta^J))_{I=0,\,\ldots,\,n}^{\mu=1,\,\ldots,\,n},\]
where
\[\chi:=\tilde{\phi}+i(\rho+c(\mathcal{K}+\log (\rho+c)))-\sum_{I=0}^n \zeta^I\tilde{\zeta}_I-\sum_{I,\,J=0}^n \zeta^IF_{IJ}(X)\zeta^J.\]
This proves the integrability of $J_1$.
\epf

 \section{Completeness of the one-loop deformation}

\subsection{Completeness of the one-loop deformation for projective special K\"ahler manifolds with regular boundary behaviour}
In this and the next section, we prove under two different types of natural assumptions the completeness of the one-loop deformed Ferrara-Sabharwal metric $g_{FS}^c$ (see Definition \ref{defDeformedFSMetric} and Theorem \ref{globalFS}) on  $\bar N$ for  $c\ge 0$. For $c<0$  and the case of projective special K\"ahler domains, $(N'_{(4n+4,\,0)},\,g_{FS}^c)$ is known to be incomplete \cite[Rem. 9]{ACDM}.\\

\bt\label{DFScompleteness}
Let $(\bar{M},\bar{g})$ be a projective special K\"ahler manifold with regular boundary behaviour and
$(\bar{N},g^c_{FS})$ the one-loop deformed Ferrara-Sabharwal (quaternionic K\"ahler) manifold associated to $(\bar{M},\bar{g})$. Then $(\bar{N},g^c_{FS})$ is complete for all $c\ge 0$. 
\et

\begin{Ex}\label{CHEx}
The projective special K\"ahler manifold $\mathbb{C}H^n$ with quadratic holomorphic prepotential $F=\frac{i}{2}((z^0)^2-\sum_{\mu=1}^n (z^\mu)^2)$ on the 
conical affine special K\"ahler domain $M:=\{|z^0|^2>\sum_{\mu=1}^n |z^\mu|^2\}$ has regular boundary behaviour in the sense of Definition 
\ref{defregbd}. Thus Corollary \ref{1stThm} implies the completeness of the projective special K\"ahler domain $\mathbb{C}H^n$. 

We know that $(\bar{N},\,g_{FS})$ is isometric to the series of Wolf spaces
\begin{equation}\tilde{X}(n+1)=\frac{SU(n+1,\,2)}{S[U(n+1)\times U(2)]}\end{equation}
of non-compact type, see e.g.\ \cite{DV}.

\bc \label{quadrCor}
For any $n\in \mathbb{N}_0$ and $c\in \mathbb{R}^{\geq 0}$,
the deformed Ferrara-Sabharwal metric
\begin{align*} g_{FS}^c&=\frac{\rho+c}{\rho}\frac{1}{1-\|X\|^2}\Big(\sum_{\mu=1}^n dX^\mu d\bar{X}^\mu+\frac{1}{1-\|X\|^2}\big|\sum_{\mu=1}^n \bar{X}^\mu dX^\mu\big|^2\Big)\\\quad&+\frac{1}{4\rho^2}\frac{\rho+2c}{\rho+c}d\rho^2-\frac{2}{\rho}(dw_0d\bar{w}_0-\sum_{\mu=1}^ndw_\mu d\bar{w}_\mu)\nonumber\\&+\frac{\rho+c}{\rho^2}\frac{4}{1-\|X\|^2}\big|dw_0+\sum_{\mu=1}^n X^\mu dw_\mu\big|^2\nonumber\\
& +\frac{1}{4\rho^2}\frac{\rho+c}{\rho+2c}\Big(d\tilde{\phi}-4\text{Im}\big(\bar{w}_0dw_0-\sum_{\mu=1}^n\bar{w}_\mu dw_\mu\big)+\frac{2c}{1-\|X\|^2}\text{Im}\sum_{\mu=1}^n \bar{X}^\mu dX^\mu\Big)^2\nonumber
\end{align*}
with $w_0:=\frac{1}{2}(\tilde{\zeta}_0+i\zeta^0)$, $w_\mu:=\frac{1}{2}(\tilde{\zeta}_\mu-i\zeta^\mu)$, $\mu=1,\,\ldots,\,n$,  on\footnote{In the case of a projective special
K\"ahler domain $\bar{M}$ we consider $\bar{N}=\bar{M}\times\mathbb{R}^{>0}\times \mathbb{R}^{2n+3}$ as in Definition \ref{defDeformedFSMetric}, rather than its cyclic quotient
$\bar{M}\times\mathbb{R}^{>0}\times S^1_c\times \mathbb{R}^{2n+2}$ on which the metric is also defined.}
\[\bar{N}=\{(X,\,\rho,\,\tilde{\phi},\,w)\in \mathbb{C}^n\times \mathbb{R}^{>0}\times\mathbb{R}\times\mathbb{C}^{n+1}\mid \|X\|^2<1\}\]
defined by the holomorphic function 
\[F=\frac{i}{2}\left((z^0)^2-\sum_{\mu=1}^n (z^\mu)^2\right)\text{ on }M:=\left\{|z^0|^2>\sum_{\mu=1}^n |z^\mu|^2\right\}\] 
is a complete quaternionic K\"ahler metric. Furthermore $(\bar{N},\,g_{FS})$ is isometric to the symmetric space $\tilde{X}(n+1)=\frac{SU(n+1,\,2)}{S[U(n+1)\times U(2)]}$.
\ec
\end{Ex}

\begin{proof}[Proof of Theorem \ref{DFScompleteness}]
Let $\gamma\colon [0,b)\to \bar{N}$ be a smooth curve which leaves every compact subset of $\bar{N}$, $b\in (0,\infty]$. We have to show that $\gamma$ has infinite length.
By Theorem \ref{2ndThm} we know that $(\bar{M},\bar{g})$ is complete. 

\bl
For every complete Riemannian manifold $(M,g)$ and $c\ge 0$ the Riemannian manifold 
\[\left(M\times \bR^{>0}, \frac{\rho+c}{\rho}g+\frac{1}{4\rho^2}\frac{\rho+2c}{\rho+c}d\rho^2\right)\]
is complete. Here $\rho$ denotes the $\bR^{>0}$-coordinate.
\el

\pf
This follows from the estimate 
\[\frac{\rho+c}{\rho}g+\frac{1}{4\rho^2}\frac{\rho+2c}{\rho+c}d\rho^2\ge g+\frac{1}{4}(d\log \rho)^2.\]
\epf

We consider the projection 
\[\bar{N}\to \bar{M}\times \bR^{>0},\; p\mapsto (\pi(p),\rho(p)),\]
where $\pi\colon \bar{N}\to \bar{M}$ is the fibre bundle projection introduced in Proposition \ref{QKGL}. Since the metric 
\[\frac{\rho+c}{\rho}g+\frac{1}{4\rho^2}\frac{\rho+2c}{\rho+c}d\rho^2\]
on the base $\bar{M}\times \bR^{>0}$ is complete by the previous lemma, the projection $\bar{\gamma}$ of $\gamma$ to $\bar{M}\times\bR^{>0}$ 
either stays in a compact set or has infinite length.  In the latter case $\gamma$ has 
infinite length. So we can assume that $\bar{\gamma}$ stays in a compact set. 

Using similar arguments as in the proof of Theorem \ref{2ndThm} we can assume that $\bar{\gamma}$ has a unique accumulation point $(\bar{p}_0,
\rho_0)$. In fact, the existence of two different accumulation points implies that
$\overline{\gamma}$ and, hence, $\gamma$ have infinite length. There exists a sequence $t_i\to b$ with $\bar{\gamma}(t_i)\to (\bar{p}_0,\rho_0)\in \bar{M}\times 
\bR^{>0}$ and $\gamma(t_i)$ leaves every compact subset of $\bar{N}_\a\cong \bar{M}_\a\times \bR^{>0} \times S^1_c\times  \bR^{2n+2}$, where $\bar{M}_\a$ is a 
projective special K\"ahler domain containing $(\bar{p}_0,\rho_0)$ and $\bar{N}_\a$ is the corresponding trivial fibre bundle endowed with the one-loop
deformed Ferrara-Sabharwal metric associated to the projective special K\"ahler domain $\bar{M}_\a$. 
Note that $\pi_{\bR^{2n+2}}(\gamma(t_i))\in \bR^{2n+2}$ is unbounded. 

\bl\label{lemmalowerbound}
For $\varepsilon >0$ and sufficiently small relatively compact $\bar{M}_\a\subset \bar{M}$ we have\footnote{Here $g_{FS}$ denotes the metric on $\bar{N}_\a=\bar{M}_\a \times \bR^{>0}\times S^1_c\times \bR^{2n+2}$ induced by the metric 
$g_{FS}$ on $\bar{M}_\a \times \bR^{>0}\times \bR^{2n+3}$. Alternatively one can compare the metrics by pulling back 
$g^c_{FS}$ to the cyclic covering $\bar{M}_\a \times \bR^{>0}\times \bR^{2n+3} \ra \bar{N}_\a$.}
$g^c_{FS}\ge \delta\cdot g_{FS}$ on $\bar{N}_\a\cap\{\rho>\varepsilon\}$ for some $\delta=\delta(\a,\varepsilon)>0$. 

\el

\pf
Choose linear coordinates $(z^0,\ldots,z^n)$ for the underlying conical affine special K\"ahler domain $M_\a$ such that $g_\a$ restricted to the 
$(z^1,\ldots,z^n)$-plane is positive definite. This can always be achieved by restricting the coordinate domain. Then it follows from \eqref{g'eq}
that $\bar{g}_\a\ge \frac{k}{4}(d^c \mathcal{K})^2$ for some $k>0$. Let $\varepsilon>0$ be given. We claim 
that 
\[g_{FS}^c\ge \frac{1}{2} \frac{k \varepsilon}{k \varepsilon+c} g_{FS}\]
on $\bar{N}_\a\cap \{\rho>\varepsilon\}$. Note first that 
\[\bar{g}+ \frac{1}{4\rho^2}\frac{\rho+2c}{\rho+c}d\rho^2\ge \frac{1}{2}\frac{k\varepsilon}{k\varepsilon +c}\left(\bar{g}
+\frac{1}{4\rho^2}d\rho^2\right).\]
Next the last two expressions in the definition of $g_{FS}^c$ can be estimated from below
\[\frac{1}{2\rho}\sum dp_a\hat{H}^{ab}dp_b+\frac{2c}{\rho^2}e^{\mathcal{K}}\left|\sum (X^Id\zt_I+F_I(X)d\z^I)\right|^2
\geq \frac{1}{2}\frac{k\epsilon}{k\epsilon+c}\frac{1}{2\rho}\sum dp_a\hat{H}^{ab}dp_b,\]
since $\frac{k\epsilon}{k\epsilon+c}\leq 1$ and $(\hat{H}^{ab})$ is positive definite.
Last setting 
\[\theta_0:=d\tilde{\phi}+\sum(\z^Id\zt_I-\zt_Id\z^I),\] 
we conclude
\begin{align*}
&\frac{c}{\rho}\bar{g}+\frac{1}{4\rho^2}\frac{\rho+c}{\rho+2c}(\theta_0+cd^c\mathcal{K})^2\nonumber\\
\ge ~~ & \frac{kc}{4\rho}(d^c\mathcal{K})^2
+ \frac{1}{4\rho^2}\underbrace{\frac{\rho+c}{\rho+2c}}_{\frac{1}{2}\leq \ldots \leq 1}\left(\underbrace{\frac{c}{k\varepsilon+c}(\theta_0+(k\varepsilon+c)d^c\mathcal{K})^2}_{\geq 0}+\frac{k\varepsilon}{k\varepsilon+c}\theta_0^2-kc\varepsilon(d^c\mathcal{K})^2\right)\nonumber\\
\geq ~~ & \frac{1}{2}\frac{k\varepsilon}{k\varepsilon+c}\frac{1}{4\rho^2}\theta_0^2+\frac{ck}{4\rho^2}(\rho-\varepsilon)(d^c\mathcal{K})^2\\
\geq ~~ &\frac{1}{2}\frac{k\varepsilon}{k\varepsilon+c}\frac{1}{4\rho^2}\left(d\tilde{\phi}+\sum(\z^Id\zt_I-\zt_Id\z^I)\right)^2,
\end{align*}
where the last inequality follows from $\rho>\varepsilon$.
Combining these three inequalities, we have shown that
\begin{equation*}
g_{FS}^c\geq \frac{1}{2}\frac{k\varepsilon}{k\varepsilon+c} g_{FS}
\end{equation*}
on $\bar{N}_\a\cap \{\rho>\varepsilon\}$.
\epf

Choose $\varepsilon>0$ such that $\rho_0\ge 2\varepsilon$. 
For the undeformed metric $g_{FS}$ on $\bar{N}_\a$ we have $g_{FS}=\bar{g}|_{\bar{M}_\a}+g_G$, where $g_G$ is a family of left invariant 
metrics on $G=\bR^{>0}\times \bR^{2n+3}$ endowed with the Lie group structure defined in \cite{CHM}. 

Since $\bar{M}_\a\subset \bar{M}$ is relatively compact, we can estimate $g_G\ge \text{const} g^0_G$ for some left invariant metric $g_G^0$ on the group fibre $G$. 
This implies that the curve $\gamma$ has infinite length, since every homogenous Riemannian metric is complete and the length of $\gamma$ can be 
estimated by the length of its projection to $G$.
\end{proof}

\subsection{Completeness of the one-loop deformation for complete projective special K\"ahler manifolds with cubic prepotential}
In this section, we prove completeness of the one-loop deformation $g^c_{FS}$ in the case of complete projective special K\"ahler manifolds in the image of the supergravity $r$-map. We will recall the definition of the latter manifolds below. They are also know as \emph{projective special K\"ahler manifolds with cubic prepotential} or 
\emph{projective very special K\"ahler manifolds}. 

\noindent
 In Section \ref{secPSRRMap}, we introduce projective special real geometry and the supergravity r-map. The latter assigns a complete projective special K\"ahler manifold to each complete projective special real manifold. In Section \ref{secCompletenessFS}, we derive a sufficient condition for the completeness of $(N'_{(4n+4,\,0)},\,g_{FS}^c)$ for $c\in\mathbb{R}^{\geq 0}$.
Recall that we construct $(N'_{(4n+4,\,0)},\,g_{FS}^c)$ from a projective special K\"ahler manifold. We prove the completeness of $(N'_{(4n+4,\,0)},\,g_{FS}^c)$ in the case that the projective special K\"ahler manifold is obtained from a complete projective special real manifold via the supergravity r-map and in the case of $\mathbb{C}H^n$. 

\noindent
As a corollary, we obtain deformations by complete quaternionic K\"ahler metrics of all known homogeneous quaternionic K\"ahler manifolds of negative scalar curvature (including symmetric spaces), except for quaternionic hyperbolic space. In the case of the series $\tilde{X}(n+1)=\frac{SU(n+1,\,2)}{S[U(n+1)\times U(2)]}$, which corresponds to the projective special K\"ahler domains $\bC H^n$ with quadratic prepotential, we already gave a simple and explicit expression for the deformed metric in Corollary \ref{quadrCor}. 

\noindent
In this chapter, we only discuss positive definite quaternionic K\"ahler metrics.

\subsubsection{Projective special real geometry and the supergravity r-map}\label{secPSRRMap}
\bd
Let $h$ be a homogeneous cubic polynomial in $n$ variables with real coefficients and let $U \subset \mathbb{R}^n\backslash\{0\}$ be an $\mathbb{R}^{>0}$-invariant domain such that $h|_U>0$ and such that $g_{\mathcal{H}}:=-\partial^2 h\big|_{\mathcal{H}}$ is a Riemannian metric on the hypersurface $\mathcal{H}:=\{x\in U\mid h(x)=1\}\subset U$. Then $(\mathcal{H},\,g_{\mathcal{H}})$ is called a \textbf{projective special real} (PSR) manifold.
\ed

\noindent
Define $\bar{M}:=\mathbb{R}^n+iU\subset\mathbb{C}^n$. We endow $\bar{M}$ with the standard complex structure and use holomorphic coordinates $(X^\mu=y^\mu+ix^\mu)_{\mu=1,\,\ldots,\,n}\in \mathbb{R}^n+iU$. We define a K\"ahler metric
\begin{align*}
\bar{g}&= 2\sum_{\mu,\,\nu=1}^ng_{\mu\bar{\nu}}dX^\mu d{\bar{X}}^\nu:= \sum_{\mu,\,\nu=1}^n\frac{\partial^2 \mathcal{K}}{\partial X^\mu \partial {\bar{X}}^\nu}dX^\mu d\bar{X}^\nu\\
&=\frac{1}{2}\sum_{\mu,\,\nu=1}^n\frac{\partial^2 \mathcal{K}}{\partial X^\mu \partial {\bar{X}}^\nu}(dX^\mu \otimes d{\bar{X}}^\nu+  d{\bar{X}}^\nu \otimes dX^\mu)\end{align*}
on $\bar{M}$ with K\"ahler potential
\be \mathcal{K}(X,\,\bar{X}):=-\text{log}\,8h(x)=-\text{log}\, h\left(i(\bar{X}-X)\right).\ee

\bd
The correspondence $(\mathcal{H},\,g_{\mathcal{H}})\mapsto (\bar{M},\,\bar{g})$ is called the \textbf{supergravity r-map}.
\ed

\br{\rmfamily\normalfont
With $\frac{\partial}{\partial X^\mu}=\frac{1}{2}\left(\frac{\partial}{\partial y^\mu}-i\frac{\partial}{\partial x^\mu}\right)$, we have
\begin{align}
2\bar{g}\left(\frac{\partial}{\partial X^\mu},\,\frac{\partial}{\partial {\bar{X}}^\nu}\right)&=2g_{\mu\bar{\nu}}=\frac{\partial^2 \mathcal{K}(X,\,\bar{X})}{\partial X^\mu \partial {\bar{X}}^\nu}=:\mathcal{K}_{\mu\bar{\nu}}\nonumber\\
&=-\frac{1}{4}\frac{\partial^2 \text{log}\,h(x)}{\partial x^\mu \partial x^\nu}=-\frac{h_{\mu\nu}(x)}{4h(x)}+\frac{h_\mu(x)h_{\nu}(x)}{4h^2(x)},\label{rMapMetric1}\end{align}
where $h_\mu(x):=\frac{\partial h(x)}{\partial x^\mu}$, $h_{\mu\nu}(x):=\frac{\partial^2 h(x)}{\partial x^\mu\partial x^\nu}$, etc.,
for $\mu,\,\nu=1,\,\ldots,\,n$.

\noindent
The inverse  $(\mathcal{K}^{\bar{\nu}\lambda})_{\nu,\,\lambda=1,\,\ldots,\,n}$ of 
$(\mathcal{K}_{\mu\bar{\nu}})_{\mu,\,\nu=1,\,\ldots,\,n}$ 
is given by 
\be \mathcal{K}^{\bar{\nu}\lambda}=-4h(x)h^{\nu\lambda}(x)+2x^\nu x^\lambda. \label{equationPVSKInverseMetric}\ee
This can be shown using the fact that $h$ is a homogeneous polynomial of degree three:
\begin{align} &\sum_{\mu=1}^n h_\mu(x)x^\mu=3h(x),\quad \sum_{\nu=1}^n h_{\mu\nu}(x)x^\nu=2h_\mu(x),\nonumber\\
&\sum_{\rho=1}^n h_{\mu\nu\rho}(x)x^\rho=h_{\mu\nu},\quad h_{\mu\nu\rho\sigma}=0. \label{HomogeneityProperty}\end{align}
}\er

\br\label{remHolPrepotrMap}{\rmfamily\normalfont
Note that any manifold $(\bar{M},\,\bar{g})$ in the image of the supergravity r-map is a projective special K\"ahler domain. The corresponding conical affine special K\"ahler domain is the trivial $\mathbb{C}^\ast$-bundle \[M:=\{z=z^0\cdot(1,\,X)\in\mathbb{C}^{n+1}\mid z^0\in\mathbb{C}^\ast,~X\in\bar{M}=\mathbb{R}^n+iU\}\to\bar{M}\] endowed with the standard complex structure $J$ and the metric $g_M$ defined by the holomorphic function
\[F:M\to \mathbb{C},\quad F(z^0,\,\ldots,\,z^n)=\frac{h(z^1,\,\cdots,\,z^n)}{z^0}.\]
Note that in general, the flat connection\footnote{$\nabla$ is defined by $x^I=\opn{Re} z^I$ and $y_I=\opn{Re} F_I(z)$ being flat, $I=0,\,\ldots,\,n$ (see \cite{ACD}).} $\nabla$ on $M$ is not the standard one induced from $\mathbb{C}^{n+1}\approx \mathbb{R}^{2n+2}$. The homothetic vector field $\xi$ is given by \linebreak$\xi=\sum_{I=0}^n(z^I\frac{\partial}{\partial z^I}+\bar{z}^I\frac{\partial}{\partial \bar{z}^I})$. To check that $\bar{g}$ is the corresponding projective special K\"ahler metric, one uses the fact that
\be 8|z^0|^2h(x)=\sum_{I,\,J=0}^nz^IN_{IJ}(z,\,\bar{z})\bar{z}^J,\ee
where as above, $x=(\opn{Im} X^1,\,\ldots,\,\opn{Im} X^n)=(\opn{Im} \frac{z^1}{z^0},\,\ldots,\,\opn{Im} \frac{z^n}{z^0})\in U$ (see \cite{CHM}).
}\er

\bd
A K\"ahler manifold $(\bar{M},\,\bar{g})$ in the image of the supergravity r-map is called a \textbf{projective very special K\"ahler manifold}.
\ed

\noindent
Due to the following {two results, projective special real geometry constitutes a powerful tool for the construction of complete projective special K\"ahler manifolds.

\bt\label{thmCHM1}{\rm \cite{CHM}}\\
The supergravity r-map preserves completeness, i.e.\ it assigns a complete projective special K\"ahler manifold to each complete projective special real manifold.
\et
The question of completeness for a projective special real manifold $(\mathcal{H},\,g_{\mathcal{H}})$ reduces to a simple topological question for the hypersurface $\mathcal{H}\subset\mathbb{R}^n$:

\bt\label{thmCNS} {\rm \cite[Thm. 2.6.]{CNS}}\\
Let $(\mathcal{H},\,g_{\mathcal{H}})$ be a projective special real manifold of dimension $n-1$. If $\mathcal{H}\subset\mathbb{R}^n$ is closed, then $(\mathcal{H},\,g_{\mathcal{H}})$ is complete.
\et

\br{\rmfamily\normalfont
In low dimensions, it is possible to classify all complete projective special real manifolds up to linear isomorphisms of the ambient space. In the case of curves, there are exactly two examples \cite{CHM}. In the case of surfaces, there exist precisely five discrete examples and a one-parameter family \cite{CDL}.
}\er

\subsubsection{The completeness theorem}\label{secCompletenessFS}
\bd
The \textbf{q-map} is the composition of the supergravity r- and c-map. It assigns a $(4n+4)$-dimensional quaternionic K\"ahler manifold to each $(n-1)$-dimensional projective special real manifold.
\ed

\br{\rmfamily\normalfont
Except for quaternionic hyperbolic space $\mathbb{H}H^{n+1}$, all Wolf spaces of non-compact type and all known homogeneous, non-symmetric quaternionic K\"ahler manifolds (called normal quaternionic K\"ahler manifolds or Alekseevsky spaces) are in the image of the supergravity c-map. While the series $\tilde{X}(n+1)=Gr_{0,\,2}(\mathbb{C}^{n+1,\,2})$ of non-compact Wolf spaces can be obtained via the supergravity c-map from the projective special K\"ahler manifold $\mathbb{C}H^n$ (with holomorphic prepotential $F=\frac{i}{2}((z^0)^2-\sum_{\mu=1}^n (z^\mu)^2)$), which is not in the image of the supergravity r-map, all the other manifolds mentioned above are in the image of the q-map.
}\er

Below, we prove the completeness of the one-loop deformation of the Ferrara-Sabharwal metric with positive deformation parameter $c\in \mathbb{R}^{\geq 0}$ for all manifolds in the image of the q-map. 

\noindent
Due to the following result, both the supergravity c-map and the q-map preserve completeness:
\bt\label{thmCMapPreservesCompleteness} {\rm \cite{CHM}}\\
The supergravity c-map assigns a complete quaternionic K\"ahler manifold of dimension $4n+4$ to each complete projective special K\"ahler manifold of dimension $2n$.
\et

Let $(\bar{M}, \bar{g})$ be a projective special K\"ahler domain with underlying conical special K\"ahler domain $(M,g,F)$. 
As in Section \ref{secHKQKCMapNew}, we assume that $M\subset \{ z^0\neq 0\} \subset \bC^{n+1}$ and identity $\bar{M}$ with 
$M\cap \{ z^0=1\}$. Then, by restricting the tensor field $\frac{g}{f}$ to $\bar{M}\subset M$, we can write
\begin{equation}
\bar{g}=-\frac{g}{f}+(\partial\mathcal{K})(\bar{\partial}\mathcal{K})=-\frac{g}{f}+\frac{1}{4}(d\mathcal{K})^2+\frac{1}{4}(d^c\mathcal{K})^2.\label{eqgPSK2}
\end{equation}

\noindent
We consider the one-loop deformed Ferrara-Sabharwal metric (see Eq.\ \eqref{DefFSmetric2})
\begin{align}g^c_{FS}=\frac{\rho+c}{\rho}\bar{g}&+\frac{1}{4\rho^2}\frac{\rho+2c}{\rho+c}d\rho^2+\frac{1}{4\rho^2}\frac{\rho+c}{\rho+2c}(d\tilde{\phi}+\sum_{I=0}^n(\z^Id\zt_I-\zt_Id\z^I)+cd^c\mathcal{K})^2\nonumber \\
&+\frac{1}{2\rho}\sum_{a,\,b=1}^{2n+2}dp_a\hat{H}^{ab}dp_b+\frac{2c}{\rho^2}e^{\mathcal{K}}\left|\sum_{I=0}^n (X^Id\zt_I+F_I(X)d\z^I)\right|^2
\end{align}
for $c\in\mathbb{R}^{\geq 0}$ defined on $N'_{(4n+4,\,0)}=\bar{N}=\bar{M}\times \mathbb{R}^{>0}\times \mathbb{R}^{2n+3}$ endowed with global coordinates \[(X^\mu,\,\rho,\,\tilde{\phi},\,\tilde{\zeta}_I,\,\zeta^I)^{\mu=1,\,\ldots,\,n}_{I=0,\,\ldots,\,n}.\]

\bp\label{propCompletenessDeformedFS2}
If $(\bar{M},\bar{g})$ is complete and $\bar{g}\geq \frac{k}{4}(d^c\mathcal{K})^2$, for some $k\in \mathbb{R}^{>0}$, then $(\bar{N},\, g_{FS}^c)$ is complete for every $c\in\mathbb{R}^{\geq 0}$.
\ep
\pf
$(\bar{N},\, g_{FS}^0)$ is complete by Theorem \ref{thmCMapPreservesCompleteness}. Since every curve on $(\bar{N},\,g_{FS}^c)$ approaching $\rho=0$ has infinite length, we can restrict to $\{\rho >\epsilon\}\subset\bar{N}$ for some $\epsilon>0$. With the same argument as in Lemma \ref{lemmalowerbound} one shows
\[g^c_{FS}\geq \frac{1}{2}\frac{k\epsilon}{k\epsilon+c}g^0_{FS}\]
using that $\bar{g}\ge \frac{k}{4}(d^c\mathcal{K})^2$.
Since $(\bar{N},\, g_{FS}^0)$ is complete, this shows that $(\bar{N},\,g_{FS}^c)$ is complete as well for $c\in\mathbb{R}^{\geq 0}$.
\epf

For quaternionic K\"ahler manifolds in the image of the q-map, the prepotential is \linebreak$F(z)=\frac{h(z^1,\,\ldots,\,z^n)}{z^0}$. 
\bl
For projective special K\"ahler manifolds in the image of the supergravity $r$-map we have
\[\bar{g}\geq \frac{1}{12}(d^c\mathcal{K})^2.\]
\el
\pf
First, we show that
\begin{equation}\tilde{g}:=-\sum_{\mu,\nu=1}^n\frac{h_{\mu\nu}(x)}{h(x)}dy^\mu dy^\nu\geq -\frac{2}{3}(d^c\mathcal{K})^2.\label{eqEstimateInLemmaCompletenessQmapDefo2}\end{equation}
Considering $\tilde{g}$ as a family of pseudo-Riemannian metrics on $\mathbb{R}^n$ depending on a parameter $x\in U$, the left hand side is positive definite on the orthogonal complement $Y^{\perp_{\tilde{g}}}$ of $Y:=\sum_{\mu=1}^n x^\mu\partial_{y^\mu}$, while the right hand side is zero, since $\tilde{g}(Y,\,\cdot)=2d^c\mathcal{K}$. In the direction of $Y$, we have
$\tilde{g}(Y,\,Y)=-6=-\frac{2}{3}(d^c\mathcal{K})^2(Y,\,Y)$.

Equation \eqref{eqEstimateInLemmaCompletenessQmapDefo2} implies
\begin{align*}
\bar{g}\geq \frac{1}{4h(x)}\sum_{\mu,\nu=1}^n \left(-h_{\mu\nu}(x)+\frac{h_\mu(x)h_\nu(x)}{h(x)}\right)dy^\mu dy^n
\geq -\frac{1}{6}(d^c\mathcal{K})^2+\frac{1}{4}(d^c\mathcal{K})^2=\frac{1}{12}(d^c\mathcal{K})^2.
\end{align*}
\epf

\noindent
This shows that the assumption of Proposition \ref{propCompletenessDeformedFS2} is fulfilled with $k=1/3$ for projective special K\"ahler manifolds in the image of the supergravity r-map and proves the following theorem. 

\bt\label{corCompleteteDefqMap2} 
Let $(\mathcal{H},\, g_{\mathcal{H}})$ be a complete projective special real manifold of dimension $n-1$ and $g_{FS}^c$, $c\in\mathbb{R}^{\geq 0}$, the one-loop deformed Ferrara-Sabharwal metric on $\bar{N}=\bar{M}\times\mathbb{R}^{>0}\times\mathbb{R}^{2n+3}$ defined by the projective special K\"ahler domain $(\bar{M},\bar{g})$ obtained from $(\mathcal{H},\, g_{\mathcal{H}})$ via the supergravity r-map. Then $(\bar{N},\,g_{FS}^c)$ is a complete quaternionic K\"ahler manifold. $(\bar{N},\,g_{FS}^0)$ is the complete quaternionic K\"ahler manifold obtained from $(\mathcal{H},\, g_{\mathcal{H}})$ via the q-map.
\et

\begin{Ex} \label{remDefG22}{\rmfamily\normalfont
For the case $n=1$ ($h=x^3$), $(\bar{N},\,g_{FS}^0)$ is isometric to the symmetric space $G_2^\ast/SO(4)$.
In this case we checked using computer algebra software that the squared pointwise norm of the Riemann tensor with respect to the metric is
\begin{align*}&\quad \sum_{i,\,j,\,k,\,l,\,\tilde{i},\,\tilde{j},\,\tilde{k},\,\tilde{l}=1}^8 R_{ijkl}g^{i\tilde{i}}g^{j\tilde{j}}g^{k\tilde{k}}g^{l\tilde{l}}R_{\tilde{i}\tilde{j}\tilde{k}\tilde{l}}\\
&=\frac{128 \left(\begin{aligned}528 c^7&+2112 c^6 \rho +3664 c^5 \rho ^2+3568 c^4 \rho ^3\\&+2110 c^3 \rho ^4+764 c^2 \rho ^5+161 c \rho ^6+17 \rho ^7\end{aligned}\right)}{3 (c+\rho ) (2 c+\rho )^6}.\end{align*}
For $c>0$, this function is non-constant, which shows that $(\bar{N},\,g_{FS}^c)$ is not locally homogeneous for $c>0$.
}
\end{Ex}

\end{document}